\documentclass[12pt,a4paper]{article}

\usepackage[utf8]{inputenc}
\usepackage[english]{babel}
\usepackage{geometry} 
\usepackage{titlesec}
\usepackage{hyperref}
\usepackage{yfonts}
\usepackage{amssymb}
\usepackage{amsmath}
\usepackage{amsthm}
\usepackage[all]{xy}
\usepackage{comment}
\usepackage{xcolor}

\newcommand{\metrica}[1]{g\left(#1\right)}

\title{Non-compact inaudibility of Naturally Reductive property}
\author{Teresa Arias-Marco and José-Manuel Fernández-Barroso}
\date{\today}

\begin{document}
\newtheorem{theorem}{Theorem}[section]
\newtheorem{proposition}[theorem]{Proposition}
\newtheorem{corollary}[theorem]{Corollary}
\newtheorem{lemma}[theorem]{Lemma}

\newtheorem*{main}{Main Theorem}
\newtheorem*{maincor1}{Main corollary 1}
\newtheorem*{maincor2}{Main corollary 2}

\theoremstyle{definition}
\newtheorem{definition}[theorem]{Definition}

\theoremstyle{definition}
\newtheorem{remark}[theorem]{Remark}

\theoremstyle{definition}
\newtheorem{example}[theorem]{Example}

\theoremstyle{definition}
\newtheorem*{notation}{Notation}

\author{Teresa Arias-Marco\footnote{ORCID: 0000-0003-0984-0367;\ email: ariasmarco@unex.es} \ and Jos\'e-Manuel Fern\'andez-Barroso\footnote{ORCID: 0000-0003-3864-9967;\ email: ferbar@unex.es}}
\date{Universidad de Extremadura, Departamento de Matemáticas, Badajoz, Spain.}

\maketitle

\begin{abstract}
 
Naturally reductive manifolds are an important class of Riemannian manifolds because they provide examples that generalize the locally symmetric ones. A property is said to be inaudible if there exists a unitary operator which intertwines the Laplace-Beltrami operator of two Riemannian manifolds such that one of them satisfies the property and the other does not.

In this paper, we study the relation between 2-step nilpotent Lie groups and the naturally reductive property to prove that this property is inaudible, using a pair of non-compact 11-dimensional generalized Heisenberg groups.
 
\end{abstract}



\section*{Introduction}

Two Riemannian manifolds $M$ and $M'$ are said to be \textit{isospectral} if there exists a unitary operator $T:L^2(M')\to L^2(M)$ which intertwines their Laplacians, that is such that $T\circ\Delta'=\Delta\circ T$.
If $M$ and $M'$ are compact, this definition is equivalent to the condition that their Laplacians have the same spectrum. This compact setting is widely studied in the literature (see \cite{AL.97,LMP.23}). In \cite{Sz.99}, Szabó constructed an operator intertwining the Laplacians of two generalized Heisenberg groups with 3-dimensional center. Similarly, the authors founded in \cite{AF.24} an operator for the 7-dimensional center case.

The well-known locally symmetric manifolds are those whose local geodesic symmetries are isometries (see \cite{H.62} for more details). The locally symmetric manifolds are also weakly symmetric, commutative, and g.o. manifolds. See, for example, the survey \cite{BTV.95} about these properties, or the comprehensive reference \cite{W.07} which provides additional background on related geometric structures. However, it is an open question whether there exists a pair of isospectral Riemannian manifolds such that one of them is locally symmetric while the other is not.

A geometric property is said to be \textit{inaudible}, or it cannot be heard, when one can find isospectral Riemannian manifolds such that one of them satisfies that property and the other does not. Gordon in \cite{G.96} noted the inaudibility of being a g.o. manifold using a pair of non-compact isospectral 23-dimensional generalized Heisenberg groups. Moreover, the authors in \cite{AF.24} used the same pair to prove that weakly symmetry and commutativity are inaudible properties on non-compact Riemannian manifolds.

Naturally reductive Riemannian manifolds $M$ are those whose geodesics in $M$ are the orbit of a one-parameter subgroup of the group of isometries, generated by a vector in the subspace $\mathfrak{m}$ of a reductive decomposition $\mathfrak{g}=\mathfrak{h}\oplus\mathfrak{m}$ of the Lie algebra $\mathfrak{g}$ of the isometry group $G$, where $\mathfrak{h}$ denotes the Lie algebra of the isotropy group $H$ of $G$. Every locally symmetric manifold is also naturally reductive. Moreover, naturally reductive manifolds are g.o. manifolds.
The classification of naturally reductive Riemannian manifolds is known up to dimension eight: in dimension three the main authors who studied this property were Tricerri and Vanhecke in \cite{TV.83}; for dimensions four and five, Kowalski and Vanhecke gave important results of their classification in \cite{KV.83} and in \cite{KV.85}, for the dimensions four and five, respectively; then, Agricola, Ferreira and Friedrich classified the six-dimensional naturally reductive spaces in \cite{AFF.15}, and more recently, Storm developed a new method in \cite{S.20} to classify naturally reductive spaces and used it to classify the seven and eight dimensional ones.

In this paper, we study the audibility of the naturally reductive property using a pair of 11-dimensional non-compact generalized Heisenberg groups. In Section \ref{sec:discussionNR-gHg} we recall the definition of generalized Heisenberg groups and we discuss when they are naturally reductive. Then, in Section \ref{sec:inaud.NR.noncompact}, we use a result of Szabó concerning isospectral non-compact generalized Heisenberg groups to set the inaudibility of being a naturally reductive manifold.

\section{Naturally reductive generalized Heisenberg groups}\label{sec:discussionNR-gHg}

Let $\mathfrak{n}=\mathfrak{v}\oplus\mathfrak{z}$, where  $\mathfrak{v}$ and $\mathfrak{z}$ are orthogonal real vector spaces with respect to an inner product $g$, and $j:\mathfrak{z}\to\mathfrak{so(v)}$ is a linear map. Then, a Lie bracket is defined on $\mathfrak{n}$ by
\begin{equation}\label{eq:relacion-j-corchete}
\metrica{[X^\mathfrak{v},Y^\mathfrak{v}]^j,Z^\mathfrak{z}}=\metrica{j_{Z^\mathfrak{z}}X^\mathfrak{v},Y^\mathfrak{v}},
\end{equation}
for $X^\mathfrak{v},Y^\mathfrak{v}\in\mathfrak{v}$ and $Z^\mathfrak{z}\in\mathfrak{z}$, such that $(\mathfrak{n},[\cdot,\cdot]^j)$ forms a 2-step nilpotent Lie algebra (i.e.  $[\mathfrak{n},\mathfrak{n}]^j\subseteq\mathfrak{z}$ and $[\mathfrak{n},\mathfrak{z}]^j=0$). We denote $(\mathfrak{n},j)$ by $\mathfrak{n}(j)$, and $(N(j),g)$ be the 2-step nilpotent Lie group whose Lie algebra is $\mathfrak{n}(j)$ with the left-invariant Riemannian metric induced by $g$, which is also denoted by $g$. The exponential map $\exp:\mathfrak{n}(j)\to N(j)$ is a diffeomorphism since $N(j)$ is simply connected and nilpotent.

When $[\mathfrak{n}(j),\mathfrak{n}(j)]^j\neq\mathfrak{z}$, the Lie group $N(j)$ is diffeomorphic to $N_1\times \mathbb{R}^k$, where $N_1=\exp(\mathfrak{v}\oplus[\mathfrak{n}(j),\mathfrak{n}(j)]^j)$ and $\mathbb{R}^k=\exp(([\mathfrak{n}(j),\mathfrak{n}(j)]^j)^\perp\cap\mathfrak{z})$. Therefore, we say that a 2-step nilpotent Lie group, $(N(j),g)$, has \textit{Euclidean factor} if it is isometric to some $(N_1,g_{|\mathfrak{n}_1\times\mathfrak{n}_1})\times\mathbb{R}^k$. Gordon proved in \cite{G.85} that a 2-step nilpotent Lie group $(N(j),g)$ has no Euclidean factor if and only if $\ker(j)=\{0\}$. 

\begin{example}
Let $\mathfrak{v}=\mathbb{R}^2$, $\mathfrak{z}=\mathbb{R}^2$ and $\{e_1,e_2,e_3,e_4\}$ an orthonormal basis of $\mathfrak{n}=\mathfrak{v}\oplus\mathfrak{z}$ with respect to an inner product $g$. For each $Z=z_3e_3+z_4e_4\in\mathfrak{z}$, consider the linear map $j:\mathfrak{z}\to\mathfrak{so}(\mathfrak{v})$ given by 
$$
j_Z=\begin{pmatrix}
 0&z_3-z_4\\
 -z_3+z_4&0
\end{pmatrix}.
$$
By \eqref{eq:relacion-j-corchete}, the only non-zero Lie bracket on $\mathfrak{n}$ is
$$
[e_1,e_2]^j=e_3-e_4.
$$
Thus, $(N(j),g)$ has Euclidean factor because $\ker(j)=\textup{span}\{(1,1)\}\neq\{0\}$. In this case, $(N(j),g)$ is isometric to the 3-dimensional Heisenberg group times a 1-dimensional Euclidean factor, $(H_3,g_{|\mathfrak{h}_3\times \mathfrak{h}_3})\times\mathbb{R}$.
\end{example}

Naturally reductive 2-step nilpotent Lie groups without Euclidean factor were characterized by Gordon in \cite{G.85}. Lauret in \cite{L.99} provided an alternative proof of this characterization using different techniques. 

\begin{theorem}[\cite{G.85,L.99}]\label{theo:caracterizacionNR}
Let $(N(j),g)$ be a 2-step nilpotent Lie group without Euclidean factor. Then, $(N(j),g)$ is naturally reductive if and only if
\begin{enumerate}
    \item $j_\mathfrak{z}=\{j_Z\}_{Z\in\mathfrak{z}}$ is a Lie subalgebra of $\mathfrak{so(v)}$.
    \item $\tau_X\in\mathfrak{so(z)}$ for any $X\in\mathfrak{z}$, where $\tau_X$ is given by $j_{X^\mathfrak{z}}j_{Y^\mathfrak{z}}-j_{Y^\mathfrak{z}}j_{X^\mathfrak{z}}=j_{\tau_{X^\mathfrak{z}}Y^\mathfrak{z}}$.
\end{enumerate}
\end{theorem}

Kaplan introduced the generalized Heisenberg groups in \cite{K.81} as special cases of 2-step nilpotent Lie groups. These manifolds are also commonly known as \textit{H-type groups} in the literature.
\begin{definition}
A \textit{generalized Heisenberg algebra} is a 2-step nilpotent Lie algebra $\mathfrak{n}(j)$ satisfying
\begin{equation}
  j_{Z^\mathfrak{z}}^2=-\|Z^\mathfrak{z}\|^2\textup{Id}_\mathfrak{v}, 
\end{equation}
for every $Z^\mathfrak{z}\in\mathfrak{z}$. The attached simply connected Lie group $(N(j),g)$ is the \textit{generalized Heisenberg group}. And such $j$ map is called a \textit{map of Heisenberg type}.
\end{definition}

The geometric information of a generalized Heisenberg group is encoded in its 2-step nilpotent Lie algebra. According to \cite{ABS.64}, the number of irreducible representations of $\mathfrak{v}$ viewed as Clifford modules together with $\dim{\mathfrak{z}}$ classifies generalized Heisenberg algebras as follows:
\begin{itemize}
\item[i)] If $\dim\mathfrak{z}\not\equiv3\mod4$, the Clifford module $Cl(\mathfrak{z})$ has a unique irreducible module $\mathfrak{v}_0$. Then, $\mathfrak{v}=(\mathfrak{v}_0)^p$ with $p\geq1$. That is, the generalized Heisenberg algebra is obtained by taking the direct sum of $p$ times the irreducible module.
\item[ii)] If $\dim\mathfrak{z}\equiv3\mod4$, the Clifford module $Cl(\mathfrak{z})$ has two non-equivalent irreducible modules, $\mathfrak{v}_1$ and $\mathfrak{v}_2$. Thus, the generalized Heisenberg algebra is obtained by taking $\mathfrak{v}=(\mathfrak{v}_1)^p\oplus(\mathfrak{v}_2)^q$, with $p\geq0,q\geq0,p+q\geq1$. We name the generalized Heisenberg algebra by $\mathfrak{n}(p,q)$ and by $N(p,q)$ its associated generalized Heisenberg group. Note that $\mathfrak{n}(p,q)$ is isomorphic to $\mathfrak{n}(q,p)$.
\end{itemize}
With this notation, if the Clifford module structure is irreducible then $\mathfrak{v}$ is said to be \textit{isotypic}. Thus, $\mathfrak{v}$ is trivially isotypic when $\dim\mathfrak{z}\not\equiv3\mod4$. If $\dim\mathfrak{z}\equiv3\mod4$, $\mathfrak{v}$ is isotypic if either $p=0$ or $q=0$.

Kaplan in \cite{K.83} classified naturally reductive generalized Heisenberg groups according to their dimension.
Tricerri and Vanhecke in \cite{TV.83} proved the same result using homogeneous structures. In both proofs, a particular map $j:\mathfrak{z}\to\mathfrak{so(v)}$ is considered given by $j_Z(X)=Z\cdot X$, with $Z\in\mathfrak{z},X\in\mathfrak{v}$, where $\cdot$ denotes the usual multiplication in $\mathfrak{v}$. Moreover, if $\mathbb{A}$ denotes the complex $\mathbb{C}$, the quaternion $\mathbb{H}$ or the Cayley (octonion) $\mathbb{O}$ numbers, $\mathfrak{v}$ is the direct sum of some copies of $\mathbb{A}$ and the center $\mathfrak{z}$ is $\mathbb{A}^*$, the non-real elements of $\mathbb{A}$. The generalized Heisenberg groups endowed with the previous $j$ maps are referred to as the Heisenberg group ($\mathfrak{z}=\mathbb{C}^*$), the quaternion analog ($\mathfrak{z}=\mathbb{H}^*$) and the Cayley analog ($\mathfrak{z}=\mathbb{O}^*$).

\begin{theorem}[\cite{K.83,TV.83}]\label{theo:NRKTV}
A generalized Heisenberg group is a naturally reductive space if and only if its center has dimension 1 (the Heisenberg group) or 3 (its quaternionic analog).
\end{theorem}

This result is previous to the characterization of naturally reductive 2-step nilmanifolds given by Gordon in \cite{G.85}.
Thus, it is necessary to clarify and specify that \textit{quaternionic analog} is equivalent to stating that $\mathfrak{v}$ is isotypic. The $j$ map used to prove Theorem \ref{theo:NRKTV} can be generalized in order to understand the isotypic and the non-isotypic generalized Heisenberg algebra at the same time. Suppose that $\dim\mathfrak{z}\equiv3\mod4$ and $\mathbb{A}$ denotes $\mathbb{H}$ or $\mathbb{O}$. Gordon introduced in \cite{G.93} the map $j:\mathfrak{z}\to\mathfrak{so(v)}$ where $\mathfrak{v}=(\mathbb{A})^p\oplus(\mathbb{A})^q$, $p,q\geq0$, $p+q\geq1$, $p,q\in\mathbb{N}$ and $\mathfrak{z}=\mathbb{A}^*$, by
$$
j_Z(X_1,\dots,X_p,X_{p+1},\dots,X_{p+q})=(Z\cdot X_1,\dots,Z\cdot X_p,X_{p+1}\cdot Z,\dots,X_{p+q}\cdot Z),
$$
where $Z\in\mathfrak{z},X_i\in \mathbb{A}, i=1,\dots p+q$, and $\cdot$ is the usual multiplication in $\mathbb{A}$. In other words, $Z\in\mathfrak{z}$ acts by the left in the first $p$ copies of $\mathbb{A}$, and it acts by the right in the last $q$ copies of $\mathbb{A}$. Note that, in the isotypic case ($p=0$ or $q=0$), this map is the same as the used by Kaplan and by Tricerri and Vanhecke. Moreover, this $j$ map is always of Heisenberg type, for every $Z\in\mathfrak{z}$ and $X\in\mathfrak{v}$, due to
$$
\begin{aligned}
j^2_Z(X_1,\dots,X_p,X_{p+1},\dots,X_{p+q})&=j_Z(Z\cdot X_1,\dots,Z\cdot X_p,X_{p+1}\cdot Z,\dots,X_{p+q}\cdot Z)\\
&=(Z^2\cdot X_1,\dots,Z^2\cdot X_p,X_{p+1}\cdot Z^2,\dots,X_{p+q}\cdot Z^2)\\
&=-\|Z\|^2\cdot(X_1,\dots,X_p,X_{p+1},\dots,X_{p+q}).
\end{aligned}
$$

In addition, $\ker(j)=\{0\}$ and the corresponding generalized Heisenberg groups with $\dim\mathfrak{z}\equiv3\mod4$ do not have Euclidean factor.

Finally, it follows that these generalized Heisenberg groups are naturally reductive if their corresponding generalized Heisenberg algebras have $\dim\mathfrak{z}=3$ and $\mathfrak{v}$ is isotypic, for instance $\mathfrak{n}=\mathfrak{n}(p,0),p\geq1$). Consider $\mathfrak{z}=\mathbb{H}^*$ and $\tau:\mathfrak{z}\times\mathfrak{z}\to\mathfrak{z}$ such that $\tau_XY= X\cdot Y-Y\cdot X$ for every orthogonal $X$ and $Y$ in $\mathfrak{z}$. Then, Theorem \ref{theo:caracterizacionNR} is satisfied, as a consequence of the properties of the quaternions and due to
$$
\begin{aligned}
j_X&j_Y(U_1,\dots,U_p)-j_Yj_X(U_1,\dots,U_p)\\
&=(X\cdot Y\cdot U_1,\dots, X\cdot Y\cdot U_p)-(Y\cdot X\cdot U_1,\dots, Y\cdot X\cdot U_p)\\
&=j_{X\cdot Y-Y\cdot X}(U_1,\dots,U_p)\\
&=j_{\tau_XY}(U_1,\dots,U_p),
\end{aligned}
$$
for every $U=(U_1,\dots,U_p)\in\mathfrak{v}$.

Now suppose that $\dim\mathfrak{z}=3$ and $\mathfrak{v}$ is not necessarily isotypic, $\mathfrak{n}=\mathfrak{n}(p,q)$, with $p,q\geq0,p+q\geq1, p,q\in\mathbb{N}$. We consider $U=U^{\mathfrak{v}_p}+U^{\mathfrak{v}_q}=(U_1,\dots,U_p,0,\dots,0)+(0,\dots,0,U_{p+1},\dots,U_{p+q})\in\mathfrak{v}$, then
$$
\begin{aligned}
j_Xj_YU-j_Yj_XU&=X\cdot Y\cdot U^{\mathfrak{v}_p}+U^{\mathfrak{v}_q}\cdot Y\cdot X-Y\cdot X\cdot U^{\mathfrak{v}_p}-U^{\mathfrak{v}_q}\cdot X\cdot Y\\
&=(X\cdot Y-Y\cdot X)\cdot U^{\mathfrak{v}_p}+U^{\mathfrak{v}_q}\cdot(Y\cdot X-X\cdot Y)\\
&=j_{X\cdot Y-Y\cdot X}U^{\mathfrak{v}_p}+j_{Y\cdot X-X\cdot Y}U^{\mathfrak{v}_q}\\
&=j_{X\cdot Y-Y\cdot X}U^{\mathfrak{v}_p}-j_{X\cdot Y-Y\cdot X}U^{\mathfrak{v}_q}\\
&=j_{X\cdot Y-Y\cdot X}(U^{\mathfrak{v}_p}-U^{\mathfrak{v}_q})
\end{aligned}
$$
which, in general, cannot be expressed in terms of $j_{\tau_XY}U$. Therefore, these generalized Heisenberg groups with 3-dimensional center and $\mathfrak{v}$ non-isotypic are not naturally reductive. Thus, the theorem proved by Kaplan in \cite{K.83} and by Tricerri and Vanhecke in \cite{TV.83} must be rewritten.
\begin{theorem}\label{theo:correccionNRHeis}
A generalized Heisenberg group is a naturally reductive space if and only if its center has dimension 1 (the Heisenberg group) or 3 with $\mathfrak{v}$ isotypic (its quaternionic analog).
\end{theorem}

\section{Non-compact inaudibility of the naturally reductivity}\label{sec:inaud.NR.noncompact}

Consider the generalized Heisenberg group $N(p,q), p,q\geq0,p+q\geq1, p,q\in\mathbb{N}$ associated with the generalized Heisenberg algebra $\mathfrak{n}(p,q)$, with 3 or 7 dimensional center. One can construct a lattice $L_{p,q}$, in $\mathfrak{n}(p,q)$, spanned by the standard basis elements. Then, $\Gamma_{p,q}=\exp(L_{p,q})$ is a cocompact discrete subgroup (i.e., it makes the quotient $N/\Gamma$ compact) of $N(p,q)$. We denote by $N^{p,q}$ the 2-step Riemannian nilmanifold $(N(p,q)/\Gamma_{p,q}, g_{p,q})$, where $g_{p,q}$ is the left-invariant Riemannian metric of $N(p,q)$. Gordon proved in \cite{G.93} the following theorem.
\begin{theorem}
If $p+q=p'+q'$, then the nilmanifolds $N^{p,q}$ and $N^{p',q'}$ are isospectral. They are locally isometric if and only if $(p',q')=\{(p,q),(q,p)\}$.
\end{theorem}

Particularly, we have the following situation
$$\xymatrix{ N(p,q)\ar[d]& N(p+q,0) \ar[d]\\ N^{p,q}\ar@{~}[r]& N^{p+q,0} }$$
where $N(p,q)$ and $N(p+q,0)$ with $p\geq0,q\geq0,p+q\geq1$, are the Riemannian covering of $N^{p,q}$ and $N^{p+q,0}$, respectively, and $\xymatrix{N^{p,q}\ar@{~}[r]& N^{p+q,0} }$ means that $N^{p,q}$ and $N^{p+q,0}$ are isospectral and not locally isometric in the compact sense. Szabó proved in \cite{Sz.99} the following result.

\begin{proposition}\label{prop:isosp-Szabo}
The generalized Heisenberg groups $N(p,q)$ and $N(p+q,0)$, $p\geq0,q\geq0,p+q\geq1$, with 3-dimensional center, are isospectral for the Laplace-Beltrami operator.
\end{proposition}

To prove it, Szabó constructed an explicit unitary operator intertwining the Laplacians of both generalized Heisenberg groups. The authors gave the same result as Szabó when $\dim\mathfrak{z}=7$, in \cite{AF.24}. Finally, we can deduce the following proposition.

\begin{theorem}
One cannot determine if a non-compact Riemannian manifold is naturally reductive using the Laplace-Beltrami operator.
\end{theorem}
\begin{proof}
Consider the generalized Heisenberg groups $N(1,1)$ and $N(2,0)$ with 3-dimensional center. By Proposition \ref{prop:isosp-Szabo}, these generalized Heisenberg groups are isospectral. Moreover, $\mathfrak{n}(2,0)$ is isotypic while $\mathfrak{n}(1,1)$ is not. Thus, using Theorem \ref{theo:correccionNRHeis}, $N(2,0)$ is naturally reductive while $N(1,1)$ is not. Therefore, the Laplace-Beltrami operator does not determine whether a non-compact Riemannian manifold is naturally reductive or not.
\end{proof}

\textbf{Authors' contributions:} All authors contributed equally to this research and in writing the paper.

\textbf{Funding:} The authors are supported by the grants GR21055 and IB18032 funded by Junta de Extremadura and Fondo Europeo de Desarrollo Regional. 
The first author is also partially supported by grant PID2019-10519GA-C22 funded by AEI/10.13039/501100011033 and by the grant GR24068 funded by Junta de Extremadura and Fondo Europeo de Desarrollo Regional.

\textbf{Conflicts of Interest:} The authors declare no conflict of interest. The founders had no role in the design of the study; in the collection, analyses, or interpretation of data; in the writing of the manuscript, or in the decision to publish the results.

\textbf{Remark:} This is a preprint of the Work accepted for publication in Siberian Mathematical Journal, \copyright, copyright 2025, Pleiades Publishing, Ltd. (\url{https://pleiades.online})


\end{document}